\documentclass[12pt,a4paper,twoside]{article}

\usepackage[applemac]{inputenc}
\usepackage[T1]{fontenc}
\usepackage[english]{babel}
\usepackage{fancyhdr}
\usepackage{indentfirst} 
\usepackage{newlfont}
\usepackage{verbatim}
\usepackage{enumerate}
\usepackage{hyperref}
\usepackage{amsthm}  
\usepackage{amssymb}  
\usepackage{amscd}
\usepackage{amsmath}
\usepackage{latexsym} 
\usepackage{amsfonts}
\usepackage{amsbsy}   
\usepackage{mathrsfs} 
\usepackage{arcs}
\usepackage{yhmath}
\usepackage{graphicx}
\usepackage{color}
\usepackage{subfigure}
\usepackage{lscape}

\theoremstyle{plain}
\newtheorem{teo}{Theorem}
\newtheorem*{teo**}{Theorem \mynumber}
\newenvironment{teo*}[1]
  {\newcommand{\mynumber}{\ref{#1}}\begin{teo**}}
  {\end{teo**}}
  
\newtheorem{prop}[teo]{Proposition}

\newtheorem{lem}[teo]{Lemma}
\newtheorem{add}[teo]{Addendum}
\theoremstyle{definition}
\newtheorem{defin}[teo]{Definition}

\theoremstyle{remark}
\newtheorem{oss}[teo]{Remark}

\begin{document}

\title{Local middle dimensional symplectic non-squeezing in the analytic setting}

\author{Lorenzo Rigolli\footnote{This work is partially supported by the DFG grant AB 360/1-1.}}
\maketitle

\begin{abstract}
We prove the following middle-dimensional non-squeezing result for analytic symplectic embeddings of domains in $\mathbb{R}^{2n}$.\\ Let $\varphi: D \hookrightarrow \mathbb{R}^{2n}$ be an analytic symplectic embedding of a domain $D \subset \mathbb{R}^{2n}$ and $P$ be a symplectic projector onto a linear $2k$-dimensional symplectic subspace $V\subset \mathbb{R}^{2n}$. Then there exists a positive function $r_0:D\rightarrow (0,+ \infty)$, bounded away from $0$ on compact subsets $K \subset D$, such that the inequality $Vol_{2k}(P\varphi (B_r(x)),\omega ^k _{0|V})\geq \pi^{k} r^{2k}$ holds for every $x \in D$ and for every $r < r_0(x)$. This claim will be deduced from an analytic middle-dimensional non-squeezing result (stated by considering paths of symplectic embeddings) whose proof will be carried on by taking advantage of a work by \'{A}lvarez Paiva and Balacheff. 
\end{abstract}

\section*{Introduction}
Let $\omega_0=\sum_{i=1}^n d x_i \wedge d y_i$ be the standard symplectic form on $\mathbb{R}^{2n}$, $B_R$ the ball of radius $R$
\begin{align*}
B_R = \{ (x_1,y_1,\ldots,x_n,y_n) \in {\mathbb{R}}^{2n} \ | \sum_{i=1} ^n x_i^2 + \sum_{i=1} ^n y_i^2  < R^2 \},
\end{align*}
and $Z_r$ the cylinder 
\begin{align*}
Z_r = \{ (x_1,y_1,\ldots,x_n,y_n) \in {\mathbb{R}}^{2n} \ | x_1 ^2 + y_1 ^2 < r^2 \}.
\end{align*}
The Gromov's non-squeezing theorem (see \cite{Gro85} or \cite{HZ94}) states that if $\varphi (B_R) \subset Z_r$, where  $\varphi$ is a symplectic (open) embedding, then $r\geq R$.
Symplectic diffeomorphisms are volume preserving, due to the fact that they preserve the multiple of the volume form ${\omega_0}^n$, but the non-squeezing theorem shows that, unlike volume preserving diffeomorphisms, they also present two-dimensional rigidity phenomena.\\
Since symplectic diffeomorphisms preserve the $2k$-form ${\omega_0}^k$ for every integer $1\leq k \leq n$, after Gromov's pioneering result one may ask if there are also middle dimensional rigidity phenomena. Some work in this direction, concerning symplectic embeddings of polydisks, has been done by Guth. In \cite{Gut08} he considers symplectic embeddings of a polydisk $\Gamma :=B_2(R_1) \times \ldots \times B_2(R_n)$ with $R_1\leq \ldots \leq R_n$ into a polydisk $\Gamma ' := B_2(R_1 ') \times \ldots \times B_2(R_n ')$ with $R_1' \leq \ldots \leq R_n '$. By the Gromov's non-squeezing theorem and the conservation of the volume under symplectic diffeomorphism, such symplectic embeddings may exist only if $R_1 \leq R_1 '$ and $R_1 \ldots	R_n \leq R_1 ' \ldots	R_n '$. On the other hand Guth proved that such embeddings exist if and only if for a certain constant $C(n)$, the inequalities $C(n) R_1 \leq R_1 '$ and $C(n) R_1 \ldots	R_n \leq R_1 ' \ldots R_n '$ hold. In particular this result exclude every middle dimensional non-squeezing phenomena in the case of polydisks embeddings.\\ 
In this paper we proceed in a different way, namely we keep the ball as domain of symplectic embeddings and in order to search for middle dimensional non-squeezing phenomena we follow the strategy pursued in \cite{AM13}.\\
First, as in \cite{EG91}, we introduce an alternative formulation of Gromov's theorem, which is that the two-dimensional shadow of the image of a radius $R$ ball in $\mathbb{R}^{2n}$ under a symplectic diffeomorphism has area at least $\pi R^2$.
More precisely the claim is that every symplectic embedding $\varphi :B_R \hookrightarrow \mathbb{R}^{2n}$ satisfies the inequality 
 \begin{align}
area( P \varphi (B_R), \omega_{0|V}) \geq \pi R^2,
\label{ineq1}
 \end{align}
where $P$ denotes the symplectic projector onto a symplectic plane $V$, i.e. the projector along the symplectic orthogonal complement of $V$.\\
This second formulation easily implies the classic one. On the other hand, if it were $area( P \varphi (B_R),\omega_{0|V}) < \pi R^2$, then, by a theorem of Moser's (see \cite{Mos65} or \cite{HZ94}), there would exist a smooth area preserving diffeomorphism\\ $\phi: P \varphi (B_R) \rightarrow B^2_r \cap V$ for some $r< R$, and then the symplectic embedding $(\phi \times id_{{V}^{\bot}}) \circ \varphi$ mapping $B_R$ into $Z_r$ would violate the classic formulation of Gromov's theorem.\\
The alternative formulation of Gromov's theorem has a natural generalization to higher dimensional shadows of a symplectic ball.\\
In other words, if $V$ is a $2k$-dimensional symplectic subspace of $\mathbb{R}^{2n}$ and $P$ is the symplectic projector onto $V$, we may ask whether it is true that 
\begin{align}
Vol_{2k} ( P \varphi(B_R), \omega_{0|V} ^k) \geq  \pi ^k R^{2k},
\label{ineq3}
\end{align}
for every symplectic embedding $\varphi: B_R \hookrightarrow \mathbb{R}^{2n}$.\\
If $k=1$ or $k=n$ the inequality holds respectively by the non-squeezing theorem and by the volume preserving property of symplectic diffeomorphisms. So we are interested only in the middle dimensional case when $2\leq k \leq n-1$.\\
If the symplectic diffeomorphism $\varphi$ is a linear map an affirmative answer to the middle dimensional non-squeezing question has been given in \cite{AM13}, nevertheless in the same paper Abbondandolo and Matveyev show that if $P$ is the symplectic projector onto a $2k$-dimensional symplectic subspace with $2\leq k \leq n-1$, then for every $\epsilon >0$ there exists an open symplectic embedding $\varphi:B_1 \hookrightarrow \mathbb{R}^{2n}$ such that $Vol_{2k} ( P \varphi(B_1), \omega_{0|V} ^k) < \epsilon$. 
Since this counterexample deforms very strongly the unitary ball, one may ask how far can the ball be deformed before the middle non-squeezing ends his validity and whether the middle dimensional non-squeezing holds locally. In \cite{AM13} the authors give two different formulations of the local question.\\
The first one asks whether, fixed a symplectic embedding $\varphi:D \hookrightarrow \mathbb{R}^{2n}$ of a domain $D\subset \mathbb{R}^{2n}$, the inequality
\begin{align}
\label{eq7}
 Vol_{2k} ( P \varphi(B_R(x)), \omega_{0|V} ^k) \geq  \pi ^k R^{2k}
\end{align}
 holds for any $x \in D$ and for $R$ positive and small enough.\\
The second formulation is the following.\\ Let us fix a path of symplectic embeddings of the unit $2n$-dimensional ball
\begin{align*}
\varphi_t :B_1 \hookrightarrow \mathbb{R}^{2n} \ \ t \in [0,1],
\end{align*}
such that $\varphi_0$ is linear (i.e. it is the restriction to $B_1$ of a linear symplectomorphism).\\
We would like to know whether there exists a positive number $t_0 \leq 1$ such that 
\begin{align}
Vol_{2k} (P \varphi_t (B_1), \omega_{0|V} ^k) \geq \pi^k, \textrm{ \ for every \ } 0\leq t < t_0.
\label{ineq5}
\end{align} 
The second formulation implies the first one by taking the path of symplectic embeddings
\begin{align}
\varphi_t (y):= 
\left\{ \begin{array}{l}
\dfrac{1}{t} \big( \varphi (x+t y) - \varphi (x) \big)  \ \ \ \textrm{ if } t \in ]0,1],\\
D \varphi (x) [y] \ \ \ \ \ \ \ \ \ \ \ \ \ \ \ \textrm{ if } t=0,
 \end{array}\right. 
\end{align} 
in fact 
\begin{align*}
Vol_{2k}(P \varphi_t (B_1(0))=Vol_{2k}(P \frac{1}{t} \varphi (B_t(x)),\omega_{0|V} ^k)=\dfrac{Vol_{2k}(P \varphi (B_t(x)),\omega_{0|V} ^k)}{t^{2k}}.
\end{align*}
In this setting Abbondandolo, Bramham, Hryniewicz and Salom\~ao \cite{ABHS15} have recently proved that if the symplectic projection is onto a $4$-dimensional symplectic subspace $V$, then both these local non-squeezing results hold (in the first formulation the diffeomorphism is required to be $C^3$). 
In this paper we address the same question but we do not impose for any assumption on the dimension of $V$, instead we require an analiticity hypothesis. First we focus on the second local formulation of the middle dimensional non-squeezing and we prove its validity under the additional assumption that the path of embeddings $t \mapsto \varphi_t$ is analytic in $t$, i.e. $\forall x \in \overline{B_1}$ the function $t\mapsto \varphi_t(x)$ is analytic.
\begin{teo}[Analytic local non-squeezing]
\label{Teor2}
Let $[0,1] \ni t \mapsto \varphi_t$ be an analytic path of symplectic embeddings $\varphi_t:\overline{B_1} \hookrightarrow \mathbb{R}^{2n}$, such that $\varphi_0$ is linear. Then the middle dimensional non-squeezing inequality
\begin{eqnarray*}
Vol_{2k} (P \varphi_t (B_1),\omega_{0|V} ^k) \geq \pi^k
\end{eqnarray*}
holds for $t$ small enough.
\end{teo} 
To prove this theorem we need some ingredients.\\
In Section \ref{Sect1} we recall some facts about contact geometry together with some results about the minimal action of a Reeb orbit in a contact manifold and we introduce Zoll contact manifolds (also known as \textit{regular contact type manifolds}), namely manifolds with the property that all Reeb orbits are periodic with the same period.\\
In Section \ref{Sect2} we prove a weaker version of the main theorem in \cite{BP14}, which says that if a constant volume deformation of the unit ball does not start tangent to all orders to a deformation by convex domains with Zoll boundaries (i.e. the deformation is \textit{not formally trivial}), then the minimal action $A_{\min}$ on the ball is strictly larger than the one of its small deformations.\\
In Section \ref{Sect4} we will see that this result implies the validity of the non-squeezing inequality \eqref{ineq5} for not formally trivial deformations of the unit ball.\\
On the other hand in case of a \textit{formally trivial deformation} we will have to proceed in a different way. Namely, using some results from Section \ref{Sect3}, we will prove that if the deformation of the ball is analytic and trivial then the function $t \mapsto Vol (P \varphi_t(B_1),\omega_{0|V} ^k)$ is analytic and has all vanishing derivatives in $t=0$. This will enable us to deduce that the function above is constant and consequently that the equality in \eqref{ineq5} holds.\\
Using Theorem \ref{Teor2} we will deduce the local non-squeezing formulation for any fixed analytic symplectic embedding. Moreover, in this latter setting we will prove that, on compact subsets of $\mathbb{R}^{2n}$, the minimal radius $R$ for which the estimate \eqref{eq7} holds is bounded away from $0$. More precisely we shall prove the following result.
\begin{teo}
\label{Teor5}
Let $\varphi: D \hookrightarrow \mathbb{R}^{2n}$ be an analytic symplectic embedding of a domain $D \subset \mathbb{R}^{2n}$. Then there exists a function $r_0:D\rightarrow (0,+ \infty)$ such that the inequality $Vol_{2k}(P \varphi (B_r(x)),\omega ^k _{0|V})\geq \pi^{k} r^{2k}$ holds, for every $x \in D$ and for every $r < r_0(x)$. Moreover $r_0$ is bounded away from $0$ on compact subsets $K \subset D$.
\end{teo}

{\textbf{Acknowledgments.}}
I would like to warmly thank Alberto Abbondandolo for all the precious help and advice he gave me concerning this paper. 
\section{Zoll contact manifolds and minimal action}
\label{Sect1}
Let us start recalling some basic facts in contact geometry.\\
A $1$-form $\alpha$ on the $(2n-1)$-dimensional manifold $M$ is a \textit{contact form} if $\alpha \wedge (d \alpha)^{n-1}$ is a volume form. In this case $(M,\alpha)$ is called \textit{contact manifold} and the volume of $M$ with respect to the volume form induced by $\alpha$ is denoted by $Vol(M,\alpha)$.\\ Moreover the contact form $\alpha$ induces a vector field $R_\alpha$ on $M$, which is called the Reeb vector field of $\alpha$ and is determined by the requirements:
\begin{eqnarray*}
i_{R_{\alpha}} d \alpha =0 \ \ \& \ \ \alpha (R_\alpha)=1.
\end{eqnarray*}

The \textit{action} $A (\gamma)$ of a periodic Reeb orbit $\gamma$ on a contact manifold $(M,\alpha)$ is defined as
\begin{align*}
A (\gamma) :=\int_{\gamma} \alpha \in \mathbb{R}.
\end{align*}
and coincides with the period of $\gamma$.
\begin{defin}
Given any contact manifold $(M,\alpha)$ with at least a closed Reeb orbit we define a function as follows
\begin{align*}
A_{\min}(M,\alpha):=  \min_{\gamma} \{  A( \gamma) \ | \  \gamma \textrm{ is a closed Reeb orbit on} \ (M,\alpha) \}.
\end{align*}
\end{defin}
Both the volume and the function $A_{\min}$ are invariant under strict contactomorphism. Indeed we have the following simple result.
\begin{prop}
Let $(M,\alpha)$ and $(N,\beta)$ be two $2n-1$ dimensional contact manifolds and $\phi:M \rightarrow N$ a strict contactomorphism (i.e. $\phi^* \beta = \alpha$).
Then
\begin{enumerate}
\item to each closed Reeb orbit $\gamma$ of $(M,\alpha)$ corresponds a closed Reeb orbit $\phi \circ \gamma$ of $(N,\beta)$,
\item $A(\gamma)= A(\phi \circ \gamma)$,
\item $A_{\min}(M,\alpha)=A_{\min}(N,\beta)$,
\item $Vol(M,\alpha)= Vol(N,\beta)$.
\end{enumerate}
\end{prop}
Another straightforward fact we will use is the following.
\begin{oss}
\label{oss1}
Let $f:S^{2n-1}\rightarrow (0,+\infty)$ be a $C^1$ function and define the set $M_f:=\{ f(x) x \ | \ x\in S^{2n-1} \} \subset \mathbb{R}^{2n}$. Then $(S^{2n-1}, f^2 \lambda_{0|S^{2n-1}})$ and $(M_f,\lambda_{0|M_f})$ are strictly contactomorphic.
\end{oss}
Indeed the radial projection $\theta:(S^{2n-1}, f^2 \lambda_{0|S^{2n-1}}) \rightarrow (M_f,\lambda_{0|M_f})$ defined by $\theta(x)=f(x)x$ is a strict contactomorphism:
\begin{align*}
& \theta^* (\lambda_{0|M_f})(x)[v]=\lambda_{0|S^{2n-1}} (\theta(x)) [ d \theta(x)[v]]=\lambda_{0|S^{2n-1}} (f(x)x)[f(x)[v]+ df(x)[v]x]=\\ 
&= f(x)\lambda_{0|S^{2n-1}} (x)[f(x)[v]]+ f(x)df(x)[v]\lambda_{0|S^{2n-1}} (x)[x]={f(x)}^2 \lambda_{0|S^{2n-1}}(x)[v].
\end{align*}
Now we introduce a special type of contact form.
\begin{defin}
A contact form on a manifold $M$ is \textit{Zoll} (or \textit{regular}) if its Reeb flow is periodic and all the Reeb orbits have the same period, and hence the same action. 
\end{defin}
For example the contact form $\lambda _{0|S^{2n-1}}$, induced on the unit sphere $S^{2n-1}$ in $\mathbb{R}^{2n}$ by the standard Liouville 1-form $\lambda_0 := \sum_{i=1}^n  x_i  d y_i$, is Zoll.\\
Later on we will consider two different kinds of deformations of a contact form: \textit{formally trivial} and \textit{not formally trivial}.
\begin{defin}
A smooth deformation $\alpha_t$, $t \in [0,t_0)$, of a contact form $\alpha_0$ is \textit{trivial} if there exists a smooth real valued function $r(t)$ and an isotopy $\phi_t$ such that $\alpha_t= r (t) \phi_t ^* \alpha_0$. 
A deformation $\alpha_t$ is \textit{formally trivial} if for every non negative $m$ there exists a trivial deformation $\alpha_t ^{(m)}$ that has order of contact $m$ with $\alpha_t$ at $t=0$. Otherwise the deformation is \textit{not formally trivial}.
\end{defin}
If instead of deformations of contact forms we choose to consider deformations of convex domains, we give the following definition. 
\begin{defin}
Consider a smooth convex domain $C_0 \subset \mathbb{R}^{2n}$ with the standard Liouville 1-form ${\lambda_0}_{\vert \partial C_0}$. 
A smooth deformation $C_t$ of $C_0$ is \textit{trivial} (resp. \textit{formally trivial}, resp. \textit{not formally trivial}) if $\theta_t ^* ( \lambda_{0|\partial C_t})$ is trivial (resp. formally trivial, resp. not formally trivial), where $\theta_t:S^{2n-1}\rightarrow \partial C_t$ is the radial projection.
\end{defin}
It is a result by Weinstein that trivial deformations of a Zoll contact form can be characterized in the following way.
\begin{prop}{\upshape{\cite{Wei74}}}
Let $\alpha_t$, $t\in [0,t_0)$, be a smooth deformation of a Zoll contact form $\alpha_0$. The deformation is trivial if and only if $\alpha_t$ is a Zoll contact form for every $t\in [0,t_0)$.
\end{prop}

In our case it will turn out that every contact deformation can be reduced to a normal form.
\begin{defin}
Let $\alpha_t = \rho_t \alpha_0$ be a smooth deformation of the Zoll contact form $\alpha_0$, where $\rho_t$ is a smooth family of positive functions on $M$, and let $m$ be a non negative integer. The deformation $\alpha_t$ is in \textit{normal form up to order $m$} if
\begin{align*}
\alpha_t = (1+ t \mu ^{(1)} + \ldots + t^m \mu ^{(m)} + t^{m+1} r_t) \alpha_0
\end{align*}
where, for $i=1,\ldots,m$, the functions $\mu ^{(i)}$ are integrals of motion for the Reeb flow of $\alpha_0$  (i.e. they are constant on the orbit of that flow) and $r_t$  is a smooth function on $M$ depending smoothly on the parameter $t$.
\end{defin}

Using a technique known as the \textit{method of Dragt and Finn} (see \cite{DF76} and \cite{Fin86}), which consist of constructing the required isotopy as composition of isotopies $\phi_t ^{(m)}$ which are flows of some particular vector fields, Balacheff and Paiva proved the following result.  
\begin{teo}{\upshape{\cite{BP14}}}
\label{Teor1}
Let $\alpha_t = \rho_t \alpha_0$ be a smooth deformation of a Zoll contact form $\alpha_0$, where $\rho_t$ is a smooth family of real valued functions on $M$ with $\rho_0=1$. Given a non negative integer $m$, there exists a contact isotopy $\phi_t ^{(m)}$ such that ${\phi_t ^{(m)}}^*\alpha _t$ is in normal form up to order $m$.
\end{teo}
\proof (Idea)
The theorem is proved by induction. The case $m=0$ follows from the Taylor expansion $\rho_t = 1+t r_t$ around $t=0$. Supposing that $\alpha_t$ is already in normal form up to order $m-1$, the point is to find a function $h_m$ in such a way that the flow $\phi_{t,{h_m}}$ of the Hamiltonian vector field $X_{h_m}$ determines an isotopy $t\mapsto \phi_{t,{h_m}}$ for which $\phi_{t,{h_m}}^* \alpha_t$ is in normal form up to order $m$.
\endproof

This proposition will be crucial in the next section:
\begin{prop}{\upshape{\cite{BP14}}}
\label{Prop3}
Let $(M,\alpha)$ be a Zoll contact manifold and let $\rho: M \rightarrow \mathbb{R}^+$ be a smooth positive function invariant under the Reeb flow of $\alpha$. Then $A_{\min} (M,\rho \alpha ) \leq A_{\min} (M,\alpha) \min \rho$.
\end{prop}
\proof
This essentially follows from the fact that if $u \in M$ is a minimum point for $\rho$, then the Reeb orbit $\gamma$ of $(M,\alpha)$ passing through $u$ is also a Reeb orbit for $(M,\rho \alpha)$. Once one checks this by a straightforward computation, we have that the action of $\gamma$ in $(M,\rho \alpha)$ is
\begin{align*}
\int_{\gamma} \rho \alpha = \min \rho \int_{\gamma} \alpha = A_{\min}  (M,\alpha) \min \rho
\end{align*}
and the proof is complete.
\endproof
We shall make use of the following classical result.
\begin{teo}{\upshape{\cite{Rab78}} \upshape{\cite{Wei78}}}
\label{Teor4}
Let $C$ be a smooth convex bounded domain of $\mathbb{R}^{2n}$.
The contact manifold $(\partial C,{\lambda_0}_{|\partial C})$ admits at least one periodic Reeb orbit.
\end{teo}
An important fact is that the function $A_{\min}$ coincides with some well known symplectic capacities such as Hofer-Zehnder and Ekeland-Hofer capacities.
\begin{teo}{\upshape{\cite{EH89}} \upshape{\cite{Vit89}}}
\label{Teor7}
Let $C$ be a smooth convex bounded domain of $\mathbb{R}^{2n}$. There exists a distinguished closed characteristic $\overline{\gamma} \subset \partial C$ such that $A_{\min}(\partial C,\lambda_0) = A( \overline{\gamma})$. Moreover, restricting to smooth convex domain of $\mathbb{R}^{2n}$, the function $A_{\min}$ is a symplectic capacity that we will denote with $c$.
\end{teo}
Due to the two theorems above, the function $A_{\min}(\partial C, {\lambda_0}_{| \partial C})$ is well defined.\\
Besides the usual proprieties of a capacity, choosing in a carefully way one among the equivalent definitions of $c$, the following result can be proved.
\begin{prop}{\upshape{\cite{AM15}}}
\label{Prop7}
Let $C \subset \mathbb{R}^{2n}$ be a smooth convex bounded domain and $P$ the symplectic projector onto a symplectic linear subspace $V \subset \mathbb{R}^{2n}$. Then $c (PC,{\omega_0}_{|V}) \geq c  (C,\omega_0)$. 
\end{prop}

\section{Deformations of $S^{2n-1}$}
\label{Sect2}
In this section we would like to get some information on how $A_{\min}$ behaves in case of a contact deformation on the unit sphere.
The results we are going to state hold in the case of an arbitrary Zoll contact manifold, but in this paper we are interested just in deformations of the standard contact form on the sphere $S^{2n-1}$, so we can simplify the proof about the Lipschitz continuity of $A_{\min}$ that relies on a result from \cite{Gin87}.
\begin{lem}
\label{Lem3}
Fix two real numbers $0<\delta<\Delta<\infty$ and consider the family $\mathcal{C}_{\delta, \Delta}$ of the convex domains $C \subset \mathbb{R}^{2n}$ which satisfy the $(\delta,\Delta)$-pinching condition $B_{\delta} \subset C  \subset B_{\Delta}$. Every symplectic capacity $c:\mathcal{C}_{\delta, \Delta} \rightarrow \mathbb{R}$ is Lipschitz continuous with respect to the Hausdorff distance.
\end{lem}
\proof
Let us take two elements $C,D \in \mathcal{C}_{\delta, \Delta}$ and let $d$ be their Hausdorff distance.\\
By assumption
\begin{eqnarray*}
\delta B=B_{\delta} \subset C,D \subset B_{\Delta}=\Delta B,
\end{eqnarray*}
hence
\begin{eqnarray*}
C \subset D + d B \subset D + \dfrac{d}{\delta} D= \left(1 + \dfrac{d}{\delta}\right) D,
\end{eqnarray*}
and by the monotonicity and conformality proprieties of symplectic capacities 
\begin{eqnarray*}
c(C) \leq \left(1+ \dfrac{d}{\delta}\right)^2 c(D) = \left(1 + 2\dfrac{d}{\delta}+\dfrac{d ^2}{\delta ^2} \right) c(D),
\end{eqnarray*}
therefore
\begin{eqnarray*}
c(C)-c(D) \leq d \left(\dfrac{2}{\delta} + \dfrac{d}{\delta ^2}\right) c(D).
\end{eqnarray*}
Because of the pinching condition we have $c(D)\leq c(\Delta B)=\Delta^2 \pi$ and $d\leq \Delta$, thus there exists a fixed real number $M>0$ such that
\begin{eqnarray*}
c(C)-c(D) \leq d M.
\end{eqnarray*} 
If $c(C) \geq c(D)$ we have $c(C)-c(D)=|c(C)-c(D)|$ and the claim follows, otherwise if $c(C) < c(D)$ we repeat the same proof switching the role of sets $C$ and $D$.
\endproof

\begin{lem}
\label{Lem2}
There exists a small open neighbourhood $U$ of zero in the Banach space $C^2 (S^{2n-1})$ such that if $f \in U$, the map $f \mapsto A_{\min} (S^{2n-1}, (1+ f) {\lambda_0}_{\vert S^{2n-1}})$ is Lipschitz continuous on $U$ with respect to the $C^0$-topology.
\end{lem}
\proof
Let us set $M_{\sqrt{1+f}}:=\{ \sqrt{1+f(x)} x \ | \ x\in S^{2n-1} \} \subset \mathbb{R}^{2n}$.\\
The map $A_{\min} (S^{2n-1}, (1+ f) \lambda_{0\vert S^{2n-1}})$ is well defined because, as observed in Remark \ref{oss1}, looking for a periodic orbit of $(S^{2n-1}, (1+ f) \lambda_{0\vert S^{2n-1}})$ is the same as looking for one of $(M_{\sqrt{1+f}}, \lambda_{0\vert M_{\sqrt{1+f}}})$, and that exists because a $C^2$-small deformation of $S^{2n-1}$ still bounds a convex domain of $\mathbb{R}^{2n}$.\\
The map $f \mapsto A_{\min}(S^{2n-1}, (1+ f) \lambda_{0\vert S^{2n-1}})$ is the composition of the maps\\ $f \mapsto M_{\sqrt{1+f}}$ and $M_{\sqrt{1+f}} \mapsto A_{\min} (M_{\sqrt{1+f}},\lambda_{0\vert M_{\sqrt{1+f}}})$, the first of which is clearly Lipschitz from the $C^0$-distance to the Hausdorff distance. So in order to prove the Lipschitz regularity result we need to show that the minimal action (which is a capacity) of a periodic orbit on a convex domain whose boundary is close to $S^{2n-1}$, is Lipschitz continuous with respect to the Hausdorff distance. But this follows from Lemma \ref{Lem3}.
\endproof

The next theorem is a weaker version, which suits in our case, of the one in \cite{BP14} which holds for every Zoll contact manifold. The proof is the same except that in our setting we do not need to use a stronger result about the Lipschitz continuity that generalizes Lemma \ref{Lem2}.
\begin{teo}
\label{Teor3}
Consider a domain $U \subset \mathbb{R}^{2n}$ and a smooth simple curve $[0,1] \ni t \mapsto y(t) \in U$ starting at $y(0)=y_0$. 
Let $(S^{2n-1},\mu_{t,y(t)})$, with $t \in [0,t_0)$, be a smooth constant volume deformation of the Zoll contact manifold $(S^{2n-1},\mu_{0,y_0}:=\lambda_{0\vert S^{2n-1}})$. If it is not formally trivial, the function $t\mapsto A_{\min} (S^{2n-1},\mu_{t,y(t)})$ attains a strict local maximum at $t=0$.
\end{teo}
\proof
To simplify the notation we will denote $\mu_{t,y(t)}$ by $\mu_t$ and, since $y(t)$ depends smoothly on $t$, we will consider that everything depends just on $t$. The proof is carried out in four steps.
\begin{enumerate}[1)]
\item First we consider the form $(1 + t \nu_t + t^m r_t) \mu_0$, where $m>1$ and both $\nu_t$ and $r_t$ are smooth function on $S^{2n-1}$ depending smoothly on $t$. By Lemma \ref{Lem2}, the function $f\mapsto A_{\min} (S^{2n-1},(1+ f)\mu_0)$ is Lipschitz if $f$ is in a small enough $C^2$-neighbourhood of zero in $C^{\infty} (S^{2n-1})$, then, for $t\rightarrow0$
\begin{align*}
A_{\min} (S^{2n-1},(1 + t \nu_t + t^m r_t) \mu_0)=A_{\min} (S^{2n-1},(1+ t \nu_t)\mu_0) + O (t^m).
\end{align*}
\item Let $(1+t^m \rho + t^{m+1} r_t) \mu_0$ be a deformation of $\mu_0$ and $\overline{\rho}$ the function obtained by averaging $\rho$ along the orbits of the Reeb vector field of $\mu_0$
\begin{align*}
\overline{\rho}(x):=\dfrac{1}{T}\int_0 ^T	 \rho (\varphi_t (x))dt,
\end{align*}
where $T$ is the common period of the periodic orbits of the Reeb flow $\varphi_t$.
According to the induction step of the proof of Theorem \ref{Teor1}, there exists a contact isotopy $\phi_t ^{(m)} :S^{2n-1} \rightarrow S^{2n-1}$ such that
\begin{align*} 
\phi_t ^{(m)*}(1+t^m \rho + t^{m+1} r_t) \mu_0=(1+t^m \overline{\rho} + t^{m+1} r' _t) \mu_0,
\end{align*} 
where $r'_t$ is a smooth function depending smoothly on $t$.
\item If $(1+t^m \rho + t^{m+1} r_t) \mu_0$ is a smooth constant volume deformation and $\overline{\rho}$ is not identically zero, then $A_{\min}(S^{2n-1},(1+t^m \rho + t^{m+1} r_t) \mu_0) <  A_{\min}(S^{2n-1},\mu_0)$ for $t\neq 0$ small enough.\\
To prove this claim, first note that by 2) and 1) follows
\begin{align*}
 A_{\min}(S^{2n-1},(1+t^m \rho + t^{m+1} r_t) \mu_0)=A_{\min}(S^{2n-1},(1+t^m \overline{\rho} + t^{m+1} r'_t) \mu_0)=
\end{align*}
\begin{align}
\label{eq3}
 = A_{\min}(S^{2n-1},(1+t^m \overline{\rho}) \mu_0)  + O({t}^{m+1}). \quad \qquad \qquad \qquad \qquad \qquad
\end{align}
 Since $\overline{\rho}$ is an integral of motion for the Reeb flow of $\mu_0$ and $m$ is a positive integer, we have that $(1 + t^m \overline{\rho})$ is a positive (for small $t$) integral of motion of $\mu_0$, thus Proposition \ref{Prop3} implies that
\begin{align}
\label{eq4}
A_{\min}(S^{2n-1},(1 + t^m \overline{\rho}) \mu_0)\leq  (1+t^m \min \overline{\rho}) A_{\min}(S^{2n-1},\mu_0).
\end{align}
The deformation $(1+t^m \overline{\rho} + t^{m+1} r'_t) \mu_0$ is constant volume because contact isotopies preserve the volume. By the proprieties of the exterior derivative 
\begin{align*}
& Vol (S^{2n-1}, \mu_0) =Vol (S^{2n-1}, (1+t^m \overline{\rho} + t^{m+1} r'_t) \mu_0) =\\
&= \int_{S^{2n-1}} {(1+t^m \overline{\rho} + t^{m+1} r'_t)}^{n} \mu_0 \wedge {d\mu_0}^{n-1}=\\
& = Vol (S^{2n-1}, \mu_0) + n t^m \int_{S^{2n-1}} \overline{\rho} \mu_0 \wedge \ {d\mu_0}^{n-1} + O({t}^{m+1}),
\end{align*}
and thus the integral of $\overline{\rho}$ over $S^{2n-1}$ is zero. Therefore, if in addition $t \neq 0$ and $\overline{\rho}$ is not identically zero, the extrema of $\overline{\rho}$ must have opposite signs and hence its minimum must be negative. Putting together this fact with \eqref{eq3} and \eqref{eq4}, we deduce that the function $t\mapsto A_{\min}(S^{2n-1},(1+t^m \rho + t^{m+1} r_t) \mu_0)$ attains a strict maximum at $t=0$, namely
\begin{align*}
A_{\min}(S^{2n-1},(1+t^m \rho + t^{m+1} r_t) \mu_0)<A_{\min}(S^{2n-1},\mu_0), \ \ \textrm{for} \ \ t>0.
\end{align*}

\item We are finally ready to prove the theorem. Let us consider a constant volume deformation $\mu_t$ of the Zoll contact form $\mu_0$. By Gray's stability theorem we can assume that the contact deformation has the form\\ $\mu_t= \rho_t \mu_0$. Expanding $\rho_t$ around $t=0$, we obtain
\begin{align*}
\mu_t = (1+ t \rho_{(1)} + t^2 r_t) \mu_0,
\end{align*} 
where $\rho_{(1)} = {\dfrac{d \rho_t}{dt}}|_{t=0}$ and $r_t$ is a smooth function depending on $t$.\\
By 3), if the average $\overline{\rho_{(1)}}$ is not identically zero, then $t\mapsto A_{\min} (S^{2n-1},\mu_t)$ attains a strict maximum at $t=0$.\\
Otherwise, if $\overline{\rho_{(1)}}$ is identically zero, by 2) there exists a contact isotopy $\phi_t^{(2)}$ such that $\phi_t ^{(2)*} \mu_t = (1+t^2 r'_t) \mu_0$. Since $\phi_t^{(2)}$ is a contact isotopy, then $(1+t^2 r'_t) \mu_0$ is also a constant volume smooth deformation of $\mu_0$ and $A_{\min} (S^{2n-1},(1+t^2 r'_t) \mu_0) =A_{\min} (S^{2n-1},\mu_t)$, so we can rewrite $\mu_t = (1+ t^2 r'_t) \mu_0$ and start anew.\\
If we repeat this process an arbitrary number of times, we see that either  $t\mapsto A_{\min} (S^{2n-1},\mu_t)$ attains a strict maximum at $t=0$ or that for any positive integer $m$, there exist a contact isotopy $\phi_t ^{(m)}$ and a smooth function $\nu_t ^{(m)}$ on $S^{2n-1}$ depending smoothly on the parameter $t$, such that $\phi_t ^{(m)*} \mu_t = (1 + t^{m} \nu_t ^{(m)} ) \mu_0$. In other words, either  $t\mapsto A_{\min} (S^{2n-1},\mu_t)$ attains a strict maximum at $t=0$ or the deformation $\mu_t$ is formally trivial.
\end{enumerate}
\endproof 

\section{Analiticity of the volume of a projection} 
\label{Sect3}
Our next goal is to prove that the fixed domain formulation of the local middle dimensional non-squeezing theorem holds if we consider an analytic path of symplectic embeddings.\\
To do this we need a result, whose proof relies on calculations made in order to prove Theorem 3 of \cite{AM13}.
\begin{prop}
\label{Prop1} Let $U\ni y_0$ be a domain of $\mathbb{R}^{n}$ and $[0,1] \times U \ni (t,y) \mapsto \varphi_{t,y}$ an analytic map such that $\varphi_{t,y}$ are embeddings of the unit $n$-dimensional ball $\varphi_{t,y}:\overline{B_1} \hookrightarrow \mathbb{R}^{n}$, with $\varphi_{0,y_0}$ linear. Moreover, let $P:\mathbb{R}^n \rightarrow V$ be the orthogonal projector onto a $k$-dimensional linear subspace $V \subset \mathbb{R}^n$ and $\rho$ a constant $k$-volume form on $V$. Then the function $(t,y) \mapsto Vol_{k} (P \varphi_{t,y} (B_1),\rho)$ is analytic in a neighbourhood of $(0,y_0)$ small enough.
\end{prop} 
In the proof we will use the following lemma.
\begin{lem} Take the hypothesis of the proposition above.
The set $S_{t,y} \subset \partial B_1$ defined as 
\begin{align}
S_{t,y} := \{x \in 	\partial B_1 | P_{| T_{\varphi_{t,y} (x)} \varphi_{t,y} (\partial B_1)} \textrm{ is not surjective} \}
\label{an1}
\end{align} 
has the property that
\begin{align}
\partial P \varphi_{t,y} (B_1)= P \varphi_{t,y} (S_{t,y})
\label{an2}
\end{align}
and can be written as
\begin{align}
S_{t,y}=\{x \in \partial B_1 | F_{t,y}(x)=0\},
\end{align}
\label{an3}
where $F_{t,y}(x):=(I-P) (D \varphi_{t,y} (x)^{*})^{-1} [x]$. If $(t,y)$ is in a small enough neighbourhood $of (0,y_0)$, $S_{t,y}$ is a submanifold of $\partial B_1$ such that $S_{t,y}= \phi_{t,y}(S^{k-1})$, where $\phi_{t,y}$ is an analytic path of diffeomorphisms.
\end{lem}
\proof
First observe that ~\eqref{an2} is an immediate consequence of the definition of $S_{t,y}$.
The function  $P_{| T_{\varphi_{t,y} (x)} \varphi_{t,y} (\partial B_1)}$ is not surjective if and only if $P D \varphi_{t,y} (x)_{| T_x \partial B_1} :T_x \partial B_1 \rightarrow T_{P \varphi_{t,y} (x)} V \cong \mathbb{R}^{k}$ is not surjective.\\
This is true iff
\begin{align*}
&\exists u \in \mathbb{R}^{k}, u \neq 0, \textrm{ such that }  <P D \varphi_{t,y}(x) [\xi], u> =0
\end{align*}
$\forall \xi \in T_x \partial B_1$, i.e. $\forall \xi$ such that  $<\xi , x> =0$.\\
Since $u = P u$ and $P=P^*$
\begin{align*}
<P D \varphi_{t,y}(x) [\xi], u>=<\xi, {(P D \varphi_{t,y}(x))}^* [u]>=<\xi, {D \varphi_{t,y}(x)}^* [u]>
\end{align*}
and thus the non surjectivity holds iff 
\begin{align*}
D \varphi_{t,y} (x) ^* [u] = \lambda x, \textrm{ where } \lambda \neq 0 \textrm{ is a real number}.
\end{align*}
Equivalently
\begin{align*}
(D \varphi_{t,y} (x)^{*})^{-1} [x] \in \mathbb{R}^{k}
\end{align*}  
which is the same as
\begin{align*}
F_{t,y}(x):=(I-P) (D \varphi_{t,y} (x)^{*})^{-1} [x] =0 \in \mathbb{R}^{n-k}.
\end{align*}  
Now, consider the analytic function $G(t,y,x):=(I-P)(\varphi_{t,y} (x)^{*})^{-1}  [x]$. We have that $\varphi_{0,y_0}= D \varphi_{0,y_0}$ because $\varphi_{0,y_0}$ is linear, hence $G(0,y,z)=0$ if $z \in S_{0,y_0}$. Applying the analytic implicit function theorem we deduce that, for $(t,y)$ close to $(0,y_0)$, $S_{t,y}$ is a submanifold of $\partial B_1$ and $S_{t,y}=\phi_{t,y} ' (S_{0,y_0})$ where $\phi_{t,y} '$ is an analytic path of diffeomorphisms. There is a diffeomorphism given by  $(D \varphi_{t,y} (x)^{*})^{-1}$ between $S^{k-1}$ and $S_{0,y_0}$, therefore by composition with $\phi_{t,y} '$ we get an analytic path of diffeomorphisms $\phi_{t,y}$ such that $\phi_{t,y}(S^{k-1})= S_{t,y}$.
\endproof
Now we are ready to prove Proposition \ref{Prop1}.   
 \proof 
Take a primitive $\alpha \in \Omega^{k-1} (V)$ of the volume form $\rho \in \Omega^{k} (V)$, i.e. $d \alpha=\rho$.\\ 
As observed in the former lemma $\partial P \varphi_{t,y} (B_1)= P \varphi_{t,y} (S_{t,y}) $ and applying Stokes' theorem we get  
\begin{align*}
Vol_k (P \varphi_{t,y} (B_1),\rho) =  \int_{P \varphi_{t,y} (B_1)} d \alpha=  \int_{\partial P \varphi_{t,y} (B_1)} \alpha=
\end{align*}
\begin{align*}
= \int_{P \varphi_{t,y} (S_{t,y})} \alpha=\int_{S_{t,y}} (P \varphi_{t,y})^* \alpha = \int_{S^{k-1}} (P \varphi_{t,y} \phi_{t,y})^* \alpha.
\end{align*}
where $\phi_{t,y}: S^{k-1} \rightarrow S_{t,y}$ is the diffeomorphism introduced in the proof of the lemma above.
For $(t,y)$ close to $(0,y_0)$, the function $(t,y) \mapsto P \varphi_{t,y} \phi_{t,y}$ is analytic and this implies the analyticity of $(t,y) \mapsto \int_{S^{k-1}} (P \varphi_{t,y} \phi_{t,y})^* \alpha$.\\
In fact, we can write $\int_{S^{k-1}} (P \varphi_{t,y} \phi_{t,y})^* \alpha=\int_{S^{k-1}} a_{t,y}(x)  \nu$ where $a_{t,y}$ is analytic.
Differentiating under integral sign, from the Taylor expansion of $a_{t,y}$ we get a local series expansion of the function $(t,y) \mapsto \int_{S^{k-1}} (P \varphi_{t,y} \phi_{t,y})^* \alpha= Vol_{k} (P \varphi_{t,y} (B_1),\rho)$, which is therefore analytic.
\endproof  
\section{Local non-squeezing}
In the following $B_1$ indicates the unit ball in $\mathbb{R}^{2n}$ and $P:\mathbb{R}^{2n} \rightarrow V$ the symplectic projection onto a $2k$-dimensional symplectic linear subspace $V \subset \mathbb{R}^{2n}$.
At first, we are interested in proving the local non squeezing formulation for a path of symplectic embeddings starting from a linear one and to do so we will use the middle dimensional linear non-squeezing result.
\label{Sect4}
\begin{teo}{\upshape{\cite{AM13}, \cite{AM15}}}
\label{teor6}
Let $P$ be the symplectic projector onto a $2k$-dimensional symplectic linear subspace $V \subset \mathbb{R}^{2n}$. Then for every linear symplectic isomorphism $L: \mathbb{R}^{2n} \rightarrow \mathbb{R}^{2n}$ there holds 
\begin{align*}
Vol_{2k} ( P L (B_1), \omega^k_{0|V}) \geq \pi^k.
\end{align*}
The equality holds if and only if the linear subspace $L^{-1} V$ is $J$-invariant, where $J$ is the standard complex structure on $\mathbb{R}^{2n}$. 
\end{teo}
We complete the above result by the following:
\begin{add}
\label{add1}
The equality holds if and only if $(P L(B_1),\omega_{0|V})$ is symplectomorphic to $(B_1 \cap L^{-1} V,\omega_{0|L^{-1} V} )$.
\end{add}
The following result is useful to prove Theorem \ref{teor6} and the addendum as well.
\begin{lem}\upshape{{\cite{Fed69} (Section 1.8.1)}}\\
\label{Lem1}
Let $1\leq k \leq n$, then 
\begin{eqnarray*}
| \omega^k[u_1,\ldots,u_{2k}]| \leq |u_1 \wedge \ldots \wedge u_{2k}| \ \ \forall u_1,\ldots , u_{2k} \in \mathbb{R}^{2n}.
\end{eqnarray*}
\end{lem}
 
\proof (Addendum)
If a symplectomorphism exists, by Lemma \ref{Lem1} we have $Vol_{2k} (P L(B_1),\omega^k _{0|V}) = Vol_{2k}(B_1 \cap L^{-1} V,\omega_{0|L^{-1} V} ^k ) \leq \pi^k$. But at the same time Theorem \ref{teor6} yields $Vol_{2k} (P L(B_1),\omega^k _{0|V}) \geq \pi^k$, hence the equality holds.  On the other hand $Vol_{2k} ( P L (B_1), \omega^k_{0|V}) = \pi^k$ iff $L^{-1} V$ is $J$-invariant; and if the claim that $PL(B_1 \cap L^{-1} V)=PL(B_1)$ is true, then $(B_1 \cap L^{-1} V,\omega_{0|L^{-1} \cap V} )$ is symplectomorphic to $(P L(B_1),\omega_{0|V})=(PL(B_1\cap L^{-1} V),\omega_{0|V})$ via the linear symplectic isomorphism $L:L^{-1} V \rightarrow V$. 
To prove the claim we reduce it to the easier case in which $P$ is orthogonal. First we take an $\omega$-compatible inner product $(\cdot,\cdot) '$ on $\mathbb{R}^{2n}$ such that $P$ is orthogonal and we denote with $B_1 '$ and $J '$ the corresponding unit ball and complex structure. In particular $V$ is $J'$-invariant. Let $\psi :(V ,\omega, J ') \rightarrow (V ,\omega, J)$ be a complex and linear isomorphism. It follows that $\psi$ is an isometry from $(V, (\cdot,\cdot)' )$ to $(V, (\cdot,\cdot) )$, hence $\psi(B_1 ')= B_1$.  The image of the unit ball under a linear surjection $M$ is given by  
\begin{align*}
M(B_1)=M( B_1 \cap rank M^*).
\end{align*}
If we take $N=L\psi$, $M=PN$ and we denote with $*'$ the adjoint of a matrix with respect to $(\cdot,\cdot)'$, we get
\begin{align*}
&PL(B_1)=PL \psi (B_1 ')= PN(B_1 ') =PN(B_1 ' \cap rank (PN)^{*'})=\\ 
&=PN(B_1 ' \cap rank (N^{*'}P^{*'}))=PN(B_1 ' \cap rank (N^{*'}P))= PN(B_1 ' \cap  N^{*'}V).
\end{align*}
The identity $\psi J' = J \psi$ implies $ J' = \psi ^{-1} J \psi$ and the fact that $L^{-1}V$ is\\ $J$-invariant is equivalent to $J\psi N^{-1}V=\psi N^{-1}V$, hence 
\begin{align*}
&N^{*'} V=N^{*'} J' V=N^{*'} J'N N^{-1} V=J' N^{-1} V=\\
&= \psi^{-1} J \psi N^{-1} V=\psi^{-1} \psi N^{-1}V=N^{-1}V,
\end{align*}
thus we obtain
\begin{align*}
&PL(B_1)=PN(B_1 ' \cap  N^{*'}V)=PN(B_1 ' \cap  N^{-1}V)=\\
&=PL\psi (B_1 ' \cap  \psi^{-1} L^{-1} V)=PL(B_1  \cap L^{-1} V),
\end{align*}
and the claim is proved.
\endproof

In order to gain some information about the strong formulation of the local non-squeezing inequality we study the function $t\mapsto Vol_{2k} (P \varphi_t (B_1),\omega_{0|V} ^k)$.
\begin{prop}
\label{prop1}
Consider a domain $U \subset \mathbb{R}^{2n}$ and a smooth simple curve $[0,1] \ni t \mapsto y(t) \in U$ starting at $y(0)=y_0$. 
Let $[0,1] \ni t \mapsto \varphi_{t,y(t)}$ be a smooth path of symplectic embeddings $\varphi_{t,y(t)}:\overline{B_1} \hookrightarrow \mathbb{R}^{2n}$, such that $\varphi_{0,y_0}$ is linear and $\varphi_{0,y_0} ^{-1} V$ is $J$-invariant.
The deformation of $P\varphi_{0,y_0}(B_1)$ given by $P\varphi_{t,y(t)}(B_1)$ can be either formally or not formally trivial:
\begin{itemize}
\item if the deformation is formally trivial, then every order $m \in \mathbb{Z}^+$ derivative of $t\mapsto Vol_{2k} (P \varphi_{t,y(t)} (B_1),\omega_{0|V} ^k)$ vanishes in $0$;
\item if the deformation is not formally trivial, then the strict middle dimensional non-squeezing inequality $Vol_{2k} (P \varphi_{t,y(t)} (B_1),\omega_{0|V} ^k) > \pi^k$ holds for $t>0$ small enough.
\end{itemize}
\end{prop}
\proof 
By the previous addendum we have that $\psi:=\varphi_{{0,y_0}|\varphi_{0,y_0} ^{-1}V}$ is a linear symplectomorphism between $(B_1 \cap \varphi_{0,y_0} ^{-1}V,\omega_{{0,y_0}|\varphi_{0,y_0} ^{-1}V})$ and $(P\varphi_{0,y_0} (B_1),\omega_{{0,y_0}|V})$.
Let us call $M_{t,y(t)}:=\partial P \varphi_{t,y(t)} (B_1)$ and consider two $1$-forms: the Liouville form $\lambda_{0|\psi^{-1} M_{t,y(t)}}$ and its pullback $\mu_{t,y(t)} := {\theta_{t,y(t)}}^* (\lambda_{0|\psi^{-1}M_{t,y(t)}})$, where $\theta_{t,y(t)} : S^{2k-1} \rightarrow \psi^{-1}M_{t,y(t)}$ is the radial diffeomorphism such that ${\theta_{t,y(t)}}^{-1}(x) = \dfrac{x}{||x||}$.\\
Later we will use the capacity $c$, which is defined only for convex domains, so let us notice once for all that, for small deformations, $\varphi_{t,y(t)} (B_1)$ is still convex and that the projection of a convex domain is still convex.\\
Now we compute the relations between the volume of the deformations.\\
Using Stokes' theorem we get
\begin{align*}
&Vol_{2k-1}(\psi^{-1}M_{t,y(t)}, \lambda_{0|\psi^{-1}M_{t,y(t)}})= \int_{\psi^{-1}\partial P \varphi_{t,y(t)} (B_1)} \lambda_{0|\psi^{-1}M_{t,y(t)}} \wedge (d \lambda_{0|\psi^{-1}M_{t,y(t)}})^{k-1}=\\ 
&=\int_{ \psi^{-1}P \varphi_{t,y(t)} (B_1)} {\omega_0}^k = Vol_{2k} (\psi^{-1}P \varphi_{t,y(t)} (B_1),\omega_{0|\varphi_{0,y_0} ^{-1}V} ^k).
\end{align*}
On the other hand, since $(\psi^{-1}M_{t,y(t)}, \lambda_{0|\psi^{-1}M_{t,y(t)}})$ and $(S^{2k-1},\mu_{t,y(t)})$ are strictly contactomorphic, $Vol_{2k-1}(\psi^{-1}M_{t,y(t)}, \lambda_{0|\psi^{-1}M_{t,y(t)}})=Vol_{2k-1}(S^{2k-1},\mu_{t,y(t)})$.
So, if $\mu_{t,y(t)} ' :=\mu_{t,y(t)} \rho (t)$, where $\rho (t):= \dfrac{1}{\sqrt[k]{Vol_{2k} (\psi^{-1}P \varphi_{t,y(t)} (B_1),\omega_{0|\varphi_{0,y_0} ^{-1} V} ^k)}}$, it follows $Vol(S^{2k-1},\mu_{t,y(t)} ')= 1$ and in particular that $\mu_{t,y(t)} '$ is a constant volume deformation.\\
Observing that closed characteristics in $(S^{2k-1}, \mu_{t,y(t)} ')$ are the same as in $(S^{2k-1}, \mu_{t,y(t)})$ we can establish the relations between the minimal action of their closed Reeb orbits
\begin{align*}
&A_{\min} (S^{2k-1}, \mu_{t,y(t)} ') = \min_{\gamma} \{ A(\gamma) \ | \ \gamma \textrm{ closed characteristic in } (S^{2k-1}, \mu_{t,y(t)} ') \} = \\
 &= \min_{\gamma}  \{ \int_{\gamma} \mu_{t,y(t)} ' \ | \ \gamma \textrm{ closed characteristic in } (S^{2k-1}, \mu_{t,y(t)} ') \} = \\ 
 &= \min_{\gamma}  \{ \rho (t)\int_{\gamma} \mu_{t,y(t)} \ | \ \gamma \textrm{ closed characteristic in } (S^{2k-1}, \mu_{t,y(t)} ') \}= \\ 
 & =\min_{\gamma}  \{ \rho (t)\int_{\gamma} \mu_{t,y(t)} \ | \ \gamma \textrm{ closed characteristic in } (S^{2k-1}, \mu_{t,y(t)} ) \}= \\
 &= \rho (t) A_{\min} (S^{2k-1}, \mu_{t,y(t)}).
\end{align*}
Since $\theta_{t,y(t)}$ is a strict contactomorphism between $(\psi^{-1}M_{t,y(t)},  \lambda_{0|\psi^{-1}M_{t,y(t)}})$ and $(S^{2k-1}, \mu_{t,y(t)} )$, we also get 
\begin{align*}
A_{\min} (S^{2k-1}, \mu_{t,y(t)} ) =A_{\min}(\psi^{-1}M_{t,y(t)},  \lambda_{0|\psi^{-1}M_{t,y(t)}})=c (\psi^{-1}P \varphi_{t,y(t)} (B_1)),
\end{align*}
where $\psi$ is a symplectomorphism.\\ Thus the quantities $A_{\min}(\psi^{-1}M_{t,y(t)},  \lambda_{0|\psi^{-1}M_{t,y(t)}})$ and $c (\psi^{-1}P \varphi_{t,y(t)} (B_1))$ are equal respectively to $A_{\min}(M_{t,y(t)},  \lambda_{0|M_{t,y(t)}})$ and  $c (P \varphi_{t,y(t)} (B_1))$.
Notice that the Weinstein conjecture holds in the convex case (Theorem \ref{Teor4}), hence a closed characteristic for $(M_{t,y(t)},  \lambda_{0|M_{t,y(t)}} )$ always exists, moreover by Theorem \ref{Teor7} the quantities above are well defined.\\
Now let us take a deformation $(S^{2k-1}, \mu_{t,y(t)} ' )$ of the standard Zoll contact form $\mu_{0,y_0}=\lambda_{0|S^{2k-1}}$ on $S^{2k-1}$, that could be either formally trivial or not formally trivial.\\ 
Suppose the former to be true, which is equivalent to say that the deformation $P\varphi_{t,y(t)}(B_1)$ is formally trivial.\\
In this case, in the last part of the proof of Theorem \ref{Teor3} we deduced that for every $m \in \mathbb{Z}^+$ there is a contact isotopy $\phi_{t,y(t)}$ such that 
\begin{align*}
\phi_{t,y(t)}^* \mu_{t,y(t)} = (1+O({t}^{m})) \mu_0.
\end{align*}
The volume function is invariant by contact isotopy, so
\begin{align*}
&Vol_{2k} (\psi^{-1}P \varphi_{t,y(t)} (B_1),\omega_{0|\varphi_{0,y_0} ^{-1}V} ^k)=Vol_{2k-1} (S^{2k-1},\mu_{t,y(t)})=\\
&=\rho(t) Vol_{2k-1} (S^{2k-1},\mu_{t,y(t)} ')=\rho(t) Vol_{2k-1} ( S^{2k-1},(1+O({t}^{m}))\mu_0), \ \ \ \forall m \in \mathbb{Z}^+.
\end{align*}
By the definition of $\rho(t)$ the above equality is equivalent to
\begin{align*}
&{Vol_{2k} (\psi^{-1}P \varphi_{t,y(t)} (B_1),\omega_{0|\varphi_{0,y_0} ^{-1}V} ^k)}^{\frac{k+1}{k}}=Vol_{2k} (\psi^{-1}P \varphi_{t,y(t)} (B_1),\omega_{0|\varphi_{0,y_0} ^{-1}V} ^k) \rho(t)=\\ 
&= Vol_{2k}( S^{2k-1},(1+O({t}^{m}))\mu_0), \ \ \ \forall m \in \mathbb{Z}^+.
\end{align*}
Therefore each of $m$-order derivatives of $Vol_{2k} (\psi^{-1}P \varphi_{t,y(t)} (B_1),\omega_{0|\varphi_{0,y_0} ^{-1}V} ^k)^{\frac{k+1}{k}}$, and hence of $Vol_{2k} (\psi^{-1}P \varphi_{t,y(t)} (B_1),\omega_{0|\varphi_{0,y_0} ^{-1}V} ^k)=Vol_{2k} (P \varphi_{t,y(t)} (B_1),\omega_{0|V} ^k)$, vanishes in $0$.\\
Now we suppose that $(S^{2k-1}, \mu_{t,y(t)} ')$ (equivalently $P\varphi_{t,y(t)}(B_1)$) is not formally trivial.
By Theorem \ref{Teor3} and the previous calculations, if $t$ is small enough the following inequality holds
\begin{align*}
& 1=\dfrac{\pi}{\sqrt[k]{\pi^k}}= \dfrac{ A_{\min} (\psi^{-1}M_{0,y_0},  \lambda_{0|\psi^{-1} M_{0,y_0}})}{\sqrt[k]{\pi^k}} = \rho (0)  A_{\min} (\psi^{-1}M_{0,y_0},  \lambda_{0|\psi^{-1}M_{0,y_0}})= \\ 
&= A_{\min} (S^{2k-1}, \mu_{0,y_0} ') >  A_{\min} (S^{2k-1}, \mu_{t,y(t)} ')=   \dfrac{A_{\min} (\psi^{-1}M_{t,y(t)},  \lambda_{0|\psi^{-1}M_{t,y(t)}})}{{\sqrt[k]{Vol_{2k} (\psi^{-1}P \varphi_{t,y(t)} (B_1),\omega_{0|\varphi_{0,y_0} ^{-1} V} ^k)}}}.
\end{align*}
So, recalling that $A_{\min} (M_{t,y(t)},  \lambda_{0|M_{t,y(t)}})=A_{\min} (\psi^{-1}M_{t,y(t)},  \lambda_{0|\psi^{-1}M_{t,y(t)}})$ and\\ $Vol_{2k} (\psi^{-1}P \varphi_{t,y(t)} (B_1),\omega_{0|\varphi_{0,y_0} ^{-1}V} ^k)=Vol_{2k} (P \varphi_{t,y(t)} (B_1),\omega_{0|V} ^k)$, if we prove that \\$A_{\min} (M_{t,y(t)},  \lambda_{0|M_{t,y(t)}})\geq \pi$, then ${ Vol_{2k} (P \varphi_{t,y(t)} (B_1),\omega_{0|V} ^k)}^{-\frac{1}{k}} < \dfrac{1}{\pi}$ and the strict local non-squeezing inequality $Vol_{2k} (P \varphi_{t,y(t)} (B_1),\omega_{0|V} ^k) > \pi^k$ holds.
But from the behaviour of the capacity $c$ respect to symplectic projections (Proposition \ref{Prop7}), we deduce
\begin{align*}
A_{\min} (M_{t,y(t)},  \lambda_{0|M_{t,y(t)}}) = c (P \varphi_{t,y(t)} (B_1)) \geq c ( \varphi_{t,y(t)} (B_1))= c ( B_1) =\pi,
\end{align*}
and hence the result.
\endproof
From this result we cannot deduce the general local non-squeezing inequality \eqref{ineq5} because in the general case we cannot say much if a trivial deformation occurs. Nevertheless, if the deformation is analytic, the local non-squeezing inequality follows easily as consequence of the proposition above.

\begin{teo*}{Teor2}[Analytic local non-squeezing]
Let $[0,1] \ni t \mapsto \varphi_t$ be an analytic path of symplectic embeddings $\varphi_t:\overline{B_1}\hookrightarrow \mathbb{R}^{2n}$, such that $\varphi_0$ is linear. Then the middle dimensional non-squeezing inequality 
\begin{align*}
Vol_{2k} (P \varphi_t (B_1),\omega_{0|V} ^k) \geq \pi^k
\end{align*}
holds for $t$ small enough.
\end{teo*}
\proof
By Theorem \ref{teor6} we have that $Vol_{2k} (P \varphi_0 (B_1),\omega_{0|V} ^k) \geq \pi^k$ and the equality holds if and only if $\varphi_0 ^{-1} V$ is $J$-invariant. If the equality does not hold the theorem is trivially true by the continuity of the volume. On the other hand, if the equality holds, Theorem \ref{teor6} implies that $\varphi_0 ^{-1}V$ is $J$-invariant and thus we are under the hypothesis of Proposition \ref{prop1}.\\ Therefore, in the case of a not formally trivial deformation $P \varphi_t(B_1)$ there is nothing to prove. Otherwise, if the deformation is formally trivial, the function $t \mapsto Vol_{2k} (P \varphi_t (B_1),\omega_{0|V} ^k)$ has vanishing derivatives in $0$, but we know by Proposition \ref{Prop1} that if $t$ is small enough this function is analytic and hence constant. Thus we get $Vol_{2k} (P \varphi_t (B_1),\omega_{0|V} ^k)=Vol_{2k} (B^{2k}_1, \omega_{0|V} ^k)= \pi^k$ for $t$ small enough.
\endproof
Note that to prove the theorem it was sufficient to use Proposition \ref{prop1} in the case where the curve $y(t)$, on which the path $t\mapsto \varphi_{t,y(t)}$ depends, is a constant curve, but the same proof leads to a generalization of Theorem \ref{Teor2} to the case in which $y(t)$ is an arbitrary analytic curve.
Thanks to this remark we can say something more about the fixed symplectic embedding formulation of the local non-squeezing, but before we state a couple of lemmata.
First a result on the local structure of the zero set of an analytic function.
\begin{teo}(Lojasiewicz's Structure Theorem) \  {\upshape{\cite[Theorem 5.2.3]{KP92}}}
\label{Teo4}
Let $f(x_1,\ldots,x_n)$ be a real analytic function in a neighbourhood of a point $y=(y_1,\ldots,y_n)$ in $\mathbb{R}^n$ and assume that $x_n \mapsto f(y_1,\ldots,y_{n-1},x_n)$ is not identically zero. There exist numbers $\delta_j>0$, $j=1,\ldots n,$ and a neighbourhood $Q_n$ (where we define $Q_k:= \{(x_1,\ldots,x_k)  \ | \ |y_j - x_j| < \delta_j, \ 1\leq j \leq k \}$) such that the zero set 
\begin{align*}
Z:=\{ x\in Q_n 	 \ | \   f(x)=0 \}
\end{align*}
has a decomposition 
\begin{align*}
Z=V^{0} \cup \ldots \cup V^{n-1},
\end{align*}
where the set $V^0$ is either empty or consists of the point $y$ alone, while for $1\leq k \leq n-1$ we may write $V^k$ as a finite disjoint union $V^k= \cup_\lambda \Gamma^k_\lambda$ of $k$-dimensional subvarieties $\Gamma^k_\lambda$. Each $\Gamma^k_\lambda$ is defined by a system of $n-k$ equations:
\begin{align*}
&x_{k+1}=^\lambda \! \! \eta _{k+1}^k (x_1,\ldots, x_k),\\
&\ \ \ \ \ \ \ \ \ \ \ldots\\
&x_n=^\lambda \! \!\eta _{n}^k (x_1,\ldots, x_k),
\end{align*}
where each function $^\lambda \eta _{k+1}^k$ is real analytic on an open subset $\Omega_\lambda ^k \subseteq Q_k \subseteq \mathbb{R}^k$.
\end{teo}
\begin{lem}
Let $\varphi: D\rightarrow \mathbb{R}^{2n}$ be an analytic symplectic embedding and $x \in D$.
As long as $x + ry \in D$, the map
\begin{align*}
\varphi_{r,x}(y) :=
\left\{ \begin{array}{l}
\dfrac{1}{r} \big( \varphi (x+ r y) - \varphi (x) \big)  \ \ \ \textrm{ if } r >0,\\
D \varphi (x) [y] \ \ \ \ \ \ \ \ \ \ \ \ \ \ \ \textrm{ if } r=0,
 \end{array}\right.
\end{align*} 
is analytic.
\end{lem}
\proof
The function $\varphi(x+ry)$ is analytic in $r$ because it is a composition of analytic maps, thus the map $\dfrac{1}{r} \big( \varphi (x+ r y) - \varphi (x) \big)$ is analytic in $r>0$.  
Since $\varphi(x+ry)- \varphi(x)$ is analytic in $r=0$, we can express it as a convergent Taylor series centred in $0$. But the $0$-th coefficient of this expansion must vanish since $\varphi(x+0y)- \varphi(x)=0$, hence we can divide by $r$ and we obtain a convergent Taylor expansion for $\dfrac{1}{r} \big( \varphi (x+ r y) - \varphi (x) \big)$ in $r=0$.  
\endproof
\begin{teo*}{Teor5}
Let $\varphi: D\hookrightarrow \mathbb{R}^{2n}$ be an analytic symplectic embedding, with $D$ domain of $\mathbb{R}^{2n}$. Then there exists a function $r_0:D\rightarrow (0,+ \infty)$ such that the inequality $Vol_{2k}(P\varphi (B_r(x)),\omega ^k_{0|V})\geq r^{2k}\pi^{k}$ holds, for every $x \in D$ and for every $r < r_0(x)$. Moreover $r_0$ is bounded away from $0$ on compact subsets $K \subset D$.
\end{teo*}
\proof
Let $\varphi_{r,x}$ be the map defined in the lemma above. Observing that
\begin{align}
\label{eq6}
Vol_{2k}(P \varphi_{r,x} (B_1(0)),\omega_{0|V} ^k)=Vol_{2k}(P \frac{1}{r} \varphi (B_r(x)),\omega_{0|V} ^k)=\dfrac{Vol_{2k}(P \varphi (B_r(x)),\omega_{0|V} ^k)}{r^{2k}},
\end{align}
for every fixed $x \in D$ we can apply Theorem \ref{Teor2} to the path $r\mapsto \varphi_{r,x}$ and we deduce the first part of the theorem.\\
Now we prove the estimate on compact sets.\\
Define a function 
\begin{align*}
f(x,r):= Vol_{2k}(P \varphi_{r,x} (B_1(0)),\omega_{0|V} ^k) - \pi^k.
\end{align*}
This function is analytic in $\mathcal{D}= \{ (x,r) \in D \times [0,+ \infty) \ | \ 0 \leq r < R(x) \}$, where $R(x)>0$ is the supremum of the radii $r$ for which $f(x,r)$ is defined. To see this is enough to apply Proposition \ref{Prop1} to the analytic map $(r,x)\mapsto \varphi_{r,x}$.
Now, take an arbitrary point $x_0 \in D$. If $f(x_0,0)>0$, then by continuity there exists a small neighbourhood $B_{\epsilon_{x_0}} \times [0,r_{x_0})$ of $(x_0,0)$ in $\mathcal{D}$, on which $f$ is positive. 
On the other hand, if $f(x_0,0)=0$, we denote with $\gamma_A ^{x_0}:[0,1] \rightarrow \mathcal{D}$ a simple analytic curve such that $\gamma_A ^{x_0}(0)=(x_0,0)$. A consequence of Theorem \ref{Teor2} is that $f(\gamma_A ^{x_0} (r))$ must be non negative in a neighbourhood of $r=0$, i.e. $(x_0,0)$ is a local minimum for the restriction of $f$ to every analytic curve $\gamma_A ^{x_0}$.
From this we can deduce that $(x_0,0)$ is a minimum for $f$ in $\mathcal{D}$.
To see it, we first extend $f$ to an analytic function in a neighbourhood of $(x_0,0)$ in $\mathbb{R}^{2n+1}$. By Theorem \ref{Teo4}, there is a small ball $B_\delta(x_0,0) \subset \mathbb{R}^{2n+1}$ in which we know how the zeros are distributed, in particular $\mathcal{D} \cap (B_\delta (x_0,0) \backslash f^{-1}(0))$ has at most a finite number $N$ of different connected components $A_i \subset \mathcal{D} \cap B_\delta (x_0,0)$ such that $(x_0,0) \in \overline{A_i}$.
The set $(f^{-1}(0) \cup_{i=1}^N A_i) \cap (\mathcal{D} \cap B_\delta (x_0,0))$ contains a neighbourhood of $(x_0,0)$ in $\mathcal{D}$, hence if we prove that $f_{|A_i}>0$ for every $i\in \{1,\ldots ,N\}$, we get the desired result.
But if it were $f_{|A_i}<0$, by Theorem \ref{Teo4} we would be able to conclude that there exists an analytic curve $\gamma_A ^{x_0}$ laying in the connected component $A_i$ and this would imply that $0$ is not a minimum for $\gamma_A ^{x_0}$, hence a contradiction.
Therefore $(x_0,0)$ is a minimum for $f$ in $\mathcal{D}$ and hence there exists a small neighbourhood $B_{\epsilon_{x_0}} \times [0,r_{x_0})$ of $(x_0,0)$ in $\mathcal{D}$ on which $f$ is positive.
Now we consider an arbitrary compact set $K \subset D$. As we have just seen, to every $x_0 \in D$ we can associate two positive real numbers $r_{x_0}$ and $\epsilon_{x_0}$. The balls of radius $\epsilon_{x_0}$ centred in an arbitrary $x_0 \in K$ produce an open cover of $K$. From this cover we can extract a finite subcover of balls of radius $\epsilon_{x_i}$ and if we define $r_0$ as the minimum in the set of the corresponding $r_{x_i}$ we get the result.
\endproof

\newpage
 \thispagestyle{empty}
\section*{\huge References}
\begin{bibliography}{A}
\bibitem[AM13]{AM13}  A. Abbondandolo and R. Matveyev,  \emph{How large is the shadow of a symplectic ball?},  J. Topol. Anal. \textbf{5} (2013), 87-119.
\bibitem[ABHS15]{ABHS15}  A. Abbondandolo, B. Bramham, U. Hryniewicz and P. Salom\~ ao \emph{Sharp systolic inequalities for Reeb flows on the three sphere},  arXiv:1054.05258 [math.SG], 2015.
\bibitem[AM15]{AM15}  A. Abbondandolo and P. Majer, \emph{A non-squeezing theorem for symplectic images of the Hilbert ball}, Calc. Var. Partial Diff. Equ. (online version).
\bibitem[BP14]{BP14}  J. C. \'{A}lvarez Paiva and F. Balacheff,  \emph{Contact geometry and isosystolic inequalities}, Geom. Funct. Anal.  \textbf{24} (2014), no. 2 648-669.
\bibitem[DF76]{DF76}  A. J. Dragt and J. M. Finn  \emph{Lie series and invariant functions for analytic symplectic maps}, J. Mathematical Phys. \textbf{17} (1976), no. 12 2215-2227.
\bibitem[EG91]{EG91}  Y. Eliashberg and M. Gromov,  \emph{Convex symplectic manifolds, Several complex variables and complex geometry, Part 2}, Proc. Sympos. Pure Math. \textbf{52} part 2 (1991), 135-162.
\bibitem[EH89]{EH89}  I. Ekeland and H. Hofer, \emph{Symplectic topology and Hamiltonian dynamics}, Math. Z. \textbf{200} (1989), 355-378.
\bibitem[Fed69]{Fed69}  H. Federer \emph{Geometric mesure theory}, Springer (1969).
\bibitem[Fin86]{Fin86}  J. M. Finn \emph{Lie transforms: a perspective}, Lecture Notes in Phys. \textbf{252}, Springer, Berlin (1986), 63-86.
\bibitem[Gin87]{Gin87}  V. L. Ginzburg \emph{New generalizations of Poincar\'e's geometric theorem}, Funktsional. Anal. i Prilozehn. \textbf{21} (1987), no.2, 16-22.
\bibitem[Gro85]{Gro85}  M. Gromov \emph{Pseudo holomorphic curves in symplectic manifolds}, Invent. Math. \textbf{82} (1985), 307-347.
\bibitem[Gut08]{Gut08}  L. Guth \emph{Symplectic embeddings of polydisks}, Invent. Math. \textbf{172} (2008), 477-489.
\bibitem[HZ94]{HZ94}  H. Hofer and E. Zehnder,  \emph{Symplectic invariants and Hamiltonian dynamics}, Birkh\"auser, 1994.
\bibitem[KP92]{KP92}  S. G. Krantz and H. R. Parks \emph{A Primer on Real Analytic Functions}, Birkh\"auser, 1992.
\bibitem[Mos65]{Mos65}  J. Moser,  \emph{On the volume elements of a manifold}, Trans. Amer. Math. Soc. \textbf{120} (1965), 286-294.
\bibitem[Rab78]{Rab78}  P. Rabinowitz, \emph{Periodic solutions of Hamiltonian Systems}, Comm. Pure Math. Appl. \textbf{31} (1978), 157-184.
\bibitem[Vit89]{Vit89}  C. Viterbo, \emph{Capacit\'e symplectiques et applications}, Ast\'erisque \textbf{177-178} (1989), no. 714 S\'eminaire Bourbaki 41\'eme ann\'ee 345-362.
\bibitem[Wei74]{Wei74}  A. Weinstein,  \emph{Fourier integral operators, quantization, and the spectra of Riemannian manifolds}, G\'eom\'etrie symplectique et physique math\'ematique, \'Editions Centre Nat. Recherche Sci., Paris (1975), 289-298.
\bibitem[Wei78]{Wei78}  A. Weinstein, \emph{Periodic orbits for convex Hamiltonian systems}, Ann. Math. \textbf{108} (1978), 507-518.

\end{bibliography} 

\newpage
 \thispagestyle{empty}

\end{document}